

\input amstex
\documentstyle {amsppt}
\document
\centerline{\bf MORDELL--WEIL PROBLEM FOR CUBIC SURFACES}

\medskip

\centerline{Yu.I.Manin}

\smallskip

\centerline{\it Max--Planck--Institut f\"ur Mathematik, Bonn, Germany}

\bigskip

\centerline{\bf \S 0. Introduction}

\medskip

\nologo
\magnification=1200
\NoBlackBoxes

Let $V$ be a plane non--singular geometrically irreducible cubic curve
over a finitely generated field $k.$ The Mordell--Weil theorem
for $V$ can be restated in the following geometric form:
there is a finite subset $B\subset V(k)$ such that the whole $V(k)$
can be obtained from $B$ by drawing secants (and tangents)
through pairs of previously constructed points and consecutively
adding their new intersection points with $V.$

\smallskip

In this note I address the question of validity of this statement
for cubic surfaces. After reminding some constructions from the book
[Ma], I analyze a numerical example, and then prove a different
version of the Mordell--Weil statement for split cubic surfaces.
A shameless change of the composition law allows me to reduce this problem to
the classical theorem on the structure of abstract projective planes.
Unfortunately, the initial question, which is more natural to ask
for minimal surfaces, remains unanswered. I would like to call attention
to this problem and its calculational aspects.

\smallskip

I am  grateful to Don Zagier whose tables are quoted in \S 2,
and to M. Rovinsky and A. Skorobogatov, discussions with whom
helped me to state and prove the main theorem.

\bigskip

\centerline{\bf \S 1. A summary of known results}

\medskip

{\bf 1.1. Notation.} Let $V$ be a cubic hypersurface without
multiple components over a field $k$ in $\bold{P}^d,\ d\ge 2.$
If $x,y,z\in V(k)$ are three points (with multiplicities)
lying on a line in $\bold{P}^d$ not belonging to $V$,
we write $x=y\circ z.$ Thus $\circ$ is a (partial and multivalued)
composition law on $V(k).$ We will also consider its restriction
on subsets of $V(k),$ e.g. that of smooth points.

\smallskip

If $x\in V(k)$ is smooth, and does not lie on a hyperplane component
of $V$, the birational map $t_x:\ V\to V,\ y\mapsto x\circ y,$
is well defined. Denote by $\roman{Bir}\ V$ the full group
of birational automorphisms ov $V.$

\smallskip

The following two results summarize the properties of $\{ t_x\}$
for curves and surfaces respectively. The first one is classical,
and the second one is proved in [M].

\smallskip

\proclaim{\quad 1.2. Theorem} Let $V$ be a smooth cubic curve. Then:

a). $\roman{Bir}\ V$ is a semidirect product of the group of projective
automorphisms and the subgroup generated by $\{ t_x\ |\ x\in V(k)\}.$

b). We have identically
$$
t_x^2=(t_xt_yt_z)^2=1 \eqno{(1.1)}
$$
for all $x,y,z\in V(k).$

\smallskip

If in addition $k$ is finitely generated over a prime field, then:

c). $\roman{Bir}\ V$ is finitely generated.

d). All points of $V(k)$ can be obtained from a finite subset of them
by drawing secants and tangents and adding the intersection
points.

\endproclaim

\smallskip

\proclaim{\quad 1.3. Theorem} Let $V$ be a minimal smooth cubic surface
over a perfect non--closed field $k.$ Then:

a). $\roman{Bir}\ V$ is a semi--direct product of the group of
projective automorphisms and the subgroup generated by
$$
\{ t_x\ |\ x\in V(k)\}\roman{\ and\ } \{ s_{u,v}\ |\ u,v\in V(K);\ [K:k]=2;\
u,v\ \roman{are\ conjugate\ over\ }k\}
$$
where
$$
s_{u,v}:=t_ut_{u\circ v}t_v.
$$

b). We have identically
$$
t_x^2=(t_xt_{x\circ y}t_y)^2=(s_{u,v})^2=1,\ st_xs^{-1}=t_{s(x)}, \eqno{(1.2)}
$$
for all pairs $u,v$ not lying on lines in $V$, and projective
automorphisms $s.$

c). The relations (1.2) form a presentation of $\roman{Bir}\ V.$

\endproclaim

\smallskip

We remind that $V$ is called {\it minimal} if one cannot blow down
some lines of $V$ by a birational morphism defined over $k.$ The
opposite class consists of {\it split} surfaces upon which all
lines are $k$--rational.

\medskip

{\bf 1.4. Discussion.} Although the two theorems are strikingly
parallel, there is an important difference between finiteness
properties in one-- and two--dimensional cases.

\smallskip

Basically, (1.1) means only that $x+y:=e\circ (x\circ y)$ is an abelian
group law with identity $e,$ whereas the statements c) and d) of
the Theorem 1.2 additionally assert that this group is finitely
generated. Therefore, (1.1) generally is not a complete
system of relations between $\{ t_x\}$.

\smallskip

Contrariwise, since (1.2) is complete, $\roman{Bir} V$ in the
twodimensional case cannot be finitely generated if $V(k)$
is infinite. This can be proved by a direct group theoretic
argument establishing a canonical form of a word in
$\{ t_x,s_{u,v}\}$ (cf. [Ma], sections 39.8.1 and 39.8.2).

\smallskip

Therefore, if something like the statement d) of Theorem 1.2
is expected to be true for cubic surfaces, this must reflect
a deep difference between relations among $\{ t_x,s_{u,v}\}$
in $\roman{Bir}\ V$ and relations among
$\{ x\}$ in $(V(k),\circ)$. The latter are much less understood than
the former. One reason is that exceptional subvarieties
of birational automorphisms are rationally parametrized
curves in $V$ which presumably should be treated as a whole
in a reasonable finiteness statement. In fact, a typical example of such
subset is a cubic curve $C(x)$ with double point $x\in V(k)$
obtained as intersection of $V$ with tangent plane at $x.$
Now, the set $(C(x)(k)\setminus\{ x\} ,\circ )$ with a composition
law $x+y=e\circ (x\circ y)$ is isomorphic to the group
of $k$--points of a form of the multiplicative group.
Such a group is not finitely generated even for $k=\bold{Q}.$
On the other hand, in $(V(k),\circ )$ this whole set must be
considered as the domain of multivalued expression $x\circ x$,
because geometrically all its points can be obtained by
drawing tangents with $k$--rational direction to $x$. Therefore
finite generation is still conceivable.

\smallskip

This comment must also help the reader to accept the definition
of a generalized operation $\circ_{(C,p)}$ in \S 3, which
is another way to deal with the same difficulty.

\bigskip

\centerline{\bf \S 2. Minimal cubic surfaces: some numerical data}

\medskip

{\bf 2.1. The structure of data.} Let $V$ be a smooth cubic
surface over a field $k$ such that $V(k)$ is infinite.
Let $h:\ V(k)\to \bold{R}_+$ be a counting function
(i.e. for all $H>0$, the set $V_H:=\{ x\in V(k)\ |\ h(x)<H\}$
is finite). In order to find a generating subset
in $(V(k),\circ )$, one can proceed as follows.

\smallskip

A. Choose a large $H$ and compile the list of all elements of $V_H$.
Let points $x$ in it be ordered by increasing $h(x).$ We will
write $x< y$ if $x$ precedes $y$, and use the number of a point
in this list as its name.

\smallskip

B. For every $x$ and every $y<x$, calculate points $x\circ y$
and choose among them those $z=x\circ y$ for which $z<x.$
Rewrite every such relation as $x=y\circ z,\ y,z<x,$ and register it
at the same line as (coordinates and number of) $x$.
Notice that if by chance $y=z$, the last relation means
exactly that $x$ lies in the tangent section of $V$ with
double point $x$.

\smallskip

If such a relation exists for $x$, we will call $x$ {\it strongly
decomposable.}

\smallskip

If all points $x$ with sufficiently large $h(x)$ were
strongly decomposable, then the ones which are not would
form a finite generating set. This is the case for cubic
curves with height as counting function. For cubic surfaces
the tables strongly indicate that it is not the case.

\smallskip

Therefore we have to consider decompositions of length $\ge 3$,
$x=M(x_1,\dots ,x_n)$, $x_i<x,$ where $M$ is a non--associative
monomial w.r.t. $\circ $. We will call {\it weakly decomposable}
points admitting such a decomposition.

\smallskip

A direct search of such decompositions is very time--consuming
(as well as a direct search of points). One problem is that
intermediate results can have height much larger than $H$;
another is that we have no a priori bound for the length
of decomposition.

\smallskip

In the example discussed below we used simple search algorithms
allowing to list those monomials $M(x_1,\dots ,x_n)<H$ for which
there is a computation scheme representing it as an iteration
of double compositions with all intermediate results registered
in $V_H$. For example, if we have two strong decompositions
$x=y\circ z=u\circ v$ with, say, $y>z,u,v,$ then we get
a weak decomposition $y=z\circ (u\circ v).$

\medskip

{\bf 2.2. An example.} D. Zagier produced a table of all primitive
solutions of $\sum_{i=1}^4ix_i^3=0$ with $h(x):=\sum_{i=1}^4|x_i|\le
1100.$ He found 379 such points and strong decompositions of 339
among them.

\smallskip

By the search described above we found weak decompositions
of further 24 points. This left us with 16 generators for
379 points, probably too many to state a finiteness conjecture.
However, there remains a possibility that this number
will diminish if decompositions with larger intermediate
results are taken into account.

\smallskip

Here are some numerical illustrations. The first three points
$1=(1,0,1,-1),\ 2=(1,1,-1,0),\ 3=(1,-1,-1,1)$ are indecomposable.
The next 26 points are strongly decomposable, e.g.
$$
24=(1,28,-19,-18)=2\circ 2=13\circ 13=14\circ 21=5\circ 23.
$$
Points 27,\ 28, and 29 are only weakly decomposable, and
$30=(15,-37,5,29)$ stubbornly resisted decomposition.

\smallskip

One of the longest decompositions found is
$$
77=5\circ (1\circ (35\circ (2\circ ( 33\circ  ((2\circ 11)\circ
(12\circ (21\circ 70))))))).
$$

\bigskip

\centerline{\bf \S3. Birationally trivial cubic surfaces: a finiteness theorem}

\medskip

{\bf 3.1. Modified composition.} Let $V$ be a smooth cubic surface,
and $x,y\in V(k).$ Let $C\subset V$ be a curve on $V$ passing
through $x,y$, and $p:\ C\to \bold{P}^2$ an embedding
of $C$ into a projective plane such that $p(C)$ is cubic, and
$p(x)\circ p(y)$ is defined in $p(C).$ We assume that $C$ and $p$
are defined over $k.$

\smallskip

In this situation we will put
$$
x\circ_{(C,p)}y:=p^{-1}(p(x)\circ p(y)).
$$

\smallskip

{\it Example 1.} Choose $C=\roman{\ a\ plane\ section\ of\ }V$
containing $x,y$. If $p$ is the embedding of $C$ into the secant
plane, then $x\circ_{(C,p)}y=x\circ y$ in the standard notation.
Notice that the result does not depend on $C$ if $x\ne y.$
If $x=y,$ then the choice of $C$ is equivalent to the choice
of a tangent line to $V$ at $x$ so that the multivaluedness
of $\circ$ is taken care of by the introduction of this new parameter.

\smallskip

{\it Example 2.} Assume now that $V$ admits a birational morphism
$p:\ V\to \bold{P}^2$ defined over $k$ (e.g., $V$ is split). We will
choose and fix $p$ once for all.
Then any plane section $C$ of $V$ not containing one of the
blown down lines as a component is embedded by $p$ into $\bold{P}^2$
as a cubic curve. Therefore we can apply to $(C,p)$ the previous
construction. Notice that this time $x\circ_{(C,p)}y$ depends
on $C$ even if $x\ne y.$

\smallskip

The following Theorem is the main result of this note:

\smallskip

\proclaim{\quad
 3.2. Theorem} Assume that $k$ is a finitely generated field.
In the situation of Example 2,
the complement $U(k)$ to the blown down lines in $V(k)$ is finitely
generated with respect to operations $\circ_{(C,p)}.$
\endproclaim

\smallskip

{\bf Proof.} Let us start with the following auxiliary construction.
Choose a $k$--rational line $l\subset \bold{P}^2.$ Then $\Gamma :=p^{-1}(l)$
is a twisted rational cubic in $V.$ The family of all such
cubics reflects properties of that of lines: a) any two different
points $a,b$ of $U(k)$ belong to a unique $\Gamma (a,b)$;\ b) any two different
$\Gamma$'s either have one common $k$--point, or intersect
a common blown down line.

\smallskip

Define now a (partial) quaternary operation on $U(k)$:
$$
*(a,b;c,d):=\Gamma (a,b)\cap \Gamma (c,d).
$$
It is defined for a Zariski dense open subset in $U(k)^4$.

\medskip

{\it Claim 1.} If $x=*(a,b;c,d)$ is well--defined, then
there exists a plane section $C$ of $V$ such that
$$
*(a,b;c,d)=a\circ_{(C,p)}b.
$$

\smallskip

In fact, choose $C$ containing $a,b,$ and $x.$ Then $p$
maps $\Gamma (a,b)$ to a line intersecting $p(C)$ at $a,b,x$.

\smallskip

It suffices now to establish the following fact:

\medskip

{\it Claim 2.} $U(k)$ is finitely generated with respect to
$*$.

\smallskip

\smallskip

To prove this, it suffices to demonstrate that $\bold{P}^2(k)$
is finitely generated with respect to the similar
quaternary operation
$$
*(a,b;c,d):=l(a,b)\cap l(c,d)
$$
where $l(a,b)$ is the line containing $a,b.$

\smallskip

In fact, start with four points in general position in $\bold{P}^2(k).$
Introduce projective coordinates using these four poits
as basic. Generate all points starting with these four and adding
intersections of lines passing through pairs of constructed points.
Obviously, the resulting set will be an abstract projective plane
satisfying the Desargues axiom. Hence it will coincide with
$\bold{P}^2(k_0)$ where $k_0$ is the prime subfield.
Represent $k$ as $k_0(t_1,\dots ,t_n).$ Add to the initial
four points the ones with coordinates $(1:t_i:0)$ and generate
a new abstract projective plane as earlier. It will contain
$\bold{P}^1(k)$ and hence coincide with $\bold{P}^2$, by a classical
reasoning: cf. [H].

\bigskip

\centerline{\bf References}

\medskip

[H] R. Hartshorne. {\it Foundations of projective geometry.} Benjamin, 1967.

\smallskip

[M] Yu. Manin. {\it Cubic forms: algebra, geometry, arithmetic.} North Holland,
1974 and 1986.

\enddocument